\documentclass[10pt]{article}

\usepackage[dvips]{graphicx}
\usepackage{amssymb}
\usepackage{amsmath}

\usepackage{caption}
\usepackage{hyperref}
\usepackage{arcs}
\usepackage {xcolor} 
\usepackage{subfigure}
\usepackage{authblk}
\usepackage{calrsfs}

\newtheorem{Proposition}{Proposition}
\newtheorem{Lemma}{Lemma}
\newtheorem{Remark}{Remark}

\newenvironment{AMS}{\small\bf 2020 AMS subject classification: }{} 

\begin{document}

\title{Evaluating Lebesgue constants by Chebyshev polynomial meshes on cube, simplex and ball}

\author[1]{L. Białas-Cież}
\affil[1]{Jagiellonian University, Poland}
\author[2]{D.J. Kenne}
\affil[2]{Doctoral School of Exact and Natural Sciences, Jagiellonian University}
\author[3]{A. Sommariva} 
\affil[3,4]{University of Padova, Italy}
\author[4]{M. Vianello}

\date{\today}
\maketitle

\begin{abstract} 

We show that product Chebyshev polynomial meshes can be used, in a fully discrete way, to evaluate with rigorous error bounds the Lebesgue constant, i.e. the maximum of the Lebesgue function, for a class of polynomial projectors on cube, simplex and ball, including interpolation, hyperinterpolation and weighted least-squares. Several examples are presented and possible generalizations outlined. A numerical software package implementing the method is freely available online.

\end{abstract}

\vskip0.2cm
\noindent
\begin{AMS}
{\rm 65D05, 65D10, 65K05.}
\end{AMS}
\vskip0.2cm
\noindent
{\small{\bf Keywords:} Multivariate polynomial meshes, cube, simplex, ball, polynomial projectors, interpolation, least-squares, hyperinterpolation, polynomial optimization, Lebesgue constant.}

\section{Introduction}

Starting from the seminal paper by Calvi and Levenberg \cite{CL08}, the notion of {\em polynomial (admissible) mesh} has been emerging in the last years as a fundamental theoretical and computational tool in multivariate polynomial approximation. Let $K\subset \mathbb{R}^d$ (or $K\subset \mathbb{C}^d$) be a polynomially determining compact set (i.e., polynomials vanishing there vanish everywhere) and denote by  $\mathbb{P}_n=\mathbb{P}_n^d$ the space of $d$-variate polynomials of degree not exceeding $n$, and by 
\begin{equation} \label{dim}
N=N_n=dim(\mathbb{P}_n)={n+d\choose n}
\end{equation}
its dimension.

A polynomial admissible mesh of $K$ is a sequence of finite norming subsets $X_n\subset K$ such that 
\begin{equation} \label{pmesh}
\|p\|_K\leq c \|p\|_{X_n}\;,\;\;\forall p\in \mathbb{P}_n\;,
\end{equation}
with $card(X_n) = O(n^\alpha)$, $\alpha\geq d$, and $c$ a constant independent
of $n$. Here and below, $\|\cdot\|_Y$ denotes the sup-norm on a continuous or discrete compact set $Y$. The fact that $card(X_n)\geq N$ necessarily holds, each $X_n$ being $\mathbb{P}_n$-determining, since polynomials vanishing on $X_n$ vanish everywhere on $K$. Such a mesh is termed {\em optimal} when $\alpha = d$. 

To give only a flavour of the topic, we recall that polynomial meshes are invariant by affine transformations, are stable under small perturbations, and can be 
assembled by finite union, finite product and algebraic transformations, 
starting from known instances. These include cubes (boxes), simplices and balls, but also more general linear and curved polytopes, as well as convex bodies and more general compact domains satysfying Markov polynomial inequalities. On the other hand, when available optimal polynomial meshes are preferable in applications due to their low cardinality, for example in the extraction of extremal points such as Fekete-like and Leja-like points. 
We do not even attempt here to give a comprehensive survey of the already considerable literature on polynomial meshes, referring the reader e.g. to \cite{BCC12,BCC16,BCLSV11,BPV20,DP23,K11,PV13,SV15} with the references therein.

In the present paper, by exploiting the connection with the polynomial optimization methods studied for example in \cite{PV18,V18,V18-2}, we develop theoretical estimates together with a numerical algorithm to approximate the value of fitting operators uniform norms (Lebesgue constants). Indeed, computation of Lebesgue constants is a matter of optimization of Lebesgue functions, which ultimately corresponds to the computation of the maximum modulus of fitting polynomials. 
The method works 
in a fully discrete way on special polynomial meshes of product Chebyshev type on cubes, simplices and balls, 
producing approximations of the Lebesgue constant from above and below, with rigorous error bounds. 
We also point out that, by the finite union property of polynomial meshes, the method can be readily applied, with the same error bounds, to complicated geometrical objects relevant in applications, namely single polytopes (via triangulation), union of polytopes and union of balls.

Computation of Lebesgue constants is important in applications, in order to check the quality of the sampling nodes 
e.g. for interpolation. In the literature this step is typically made, for example within spectral and high-order methods for numerical PDEs, by empirical approaches, namely by evaluation of Lebesgue functions on increasingly dense discretizations. In the present paper we provide, apparently for the first time, a fully discrete method with rigorous error bounds as well as the corresponding numerical codes.

The paper is organized as follows. In Section 2 we state and prove the main theoretical results and we also outline an extended (but less accurate) approach via general polynomial meshes on multidimensional compact sets. In Section 3 we discuss some computational and implementation issues, and in Section 4 we present several examples concerning evaluation of the size of Lebesgue constants for polynomial interpolation and least-squares (including hyperinterpolation) on cube, simplex and ball. 

\section{Lebesgue constants by Chebyshev meshes}

In the sequel, the following constant will play a key role
$$
c_m=\frac{1}{\cos(\pi/(2m))}\;,\;\;m>1\;.
$$
Moreover, we shall denote by $\mathcal{C}_k$ the set of $k$ Chebyshev zeros in $(-1,1)$, 
$\cos((2j-1)\pi/(2k))$, $1\leq j\leq k$ (Chebyshev points), or the set of $k+1$ Chebyshev extrema 
in $[-1,1]$, $\cos(j\pi/k)$, $0\leq j\leq k$ (Chebyshev-Lobatto points).

\subsection{The cube}

We discuss first the case of the $d$-cube. By invertible affine transformation, the result can be immediately extended to any $d$-box, 
with invariance of the mesh constant $c$. The following result has been proved in \cite{PV18} essentially following \cite{B17}, by the notion of {\em Dubiner distance} on a compact set (which is tailored to polynomial spaces; cf. e.g. \cite{BLW08,V19}) with the references therein).

\begin{Lemma}
Let $\mathcal{C}_k$ be Chebyshev or Chebyshev-Lobatto points in $[-1,1]$. Then the sequence of product Chebyshev grids $X_n^m=(\mathcal{C}_{mn})^d$, $n=1,2,\dots$, for a fixed $m>1$, is an admissible polynomial mesh for the $d$-cube $[-1,1]^d$, with constant $c=c_m$.
\end{Lemma}

We can now state and prove a basic result on approximation of Lebesgue constants by polynomial meshes on the $d$-cube. 

\begin{Proposition}
Let $X_n^m=(\mathcal{C}_{mn})^d$ be a product Chebyshev admissible mesh of $K=[-1,1]^d$ as in Lemma 1, and $L_n:C(K)\to \mathbb{P}_n$ a linear projection operator such that 
\begin{equation} \label{Ln}
L_nf(x)=\sum_{i=1}^M{f(\xi_i)\,\varphi_i(x)}\;,
\end{equation}
where $\Xi=\{\xi_i\}\subset K$ and $\{\varphi_i\}$   is a set of generators of $\mathbb{P}_n$. Moreover, let 
$$
\lambda_n(x)=\sum_{i=1}^M{|\varphi_i(x)|}
$$
be the ``Lebesgue function''of $L_n$ and 
 $$
 \|L_n\|=\sup_{f\neq 0}\frac{\|L_nf\|_K}{\|f\|_K}=\|\lambda_n\|_K
 $$
 its ``Lebesgue constant''. 
Then, the following estimate holds for every $m>1$ 
\begin{equation} \label{lebest}
\|\lambda_n\|_{X_n^m}\leq \|L_n\|\leq 
c_m\|\lambda_n\|_{X_n^m}\;.
\end{equation} 

\end{Proposition}
\vskip0.3cm
\noindent{\em Proof.} 
First we prove that $\|L_n\|=\|\lambda_n\|_K$ for any projection operator of the form (\ref{Ln}), on a general compact set $K$. 
Indeed, inequality $\|L_nf\|_K\leq \|f\|_K \|\lambda_n\|_K$ is immediate, whereas existence of $g^\ast\in C(K)$ such that 
$g^\ast(\xi_i)=sign(\varphi_i(x^\ast))$, where $\|\lambda_n\|_K=\lambda_n(x^\ast)$ with $x^\ast \in K$, and $\|g^\ast\|_K=1$, is guaranteed by a  quite general topological result, namely the celebrated Tietze extension theorem of continuous functions from a closed subset of a normal topological space, preserving the range (cf. e.g. \cite[Ch.7, Thm.5.1]{D66}). Indeed, defining a function $g$ on the sampling set $\Xi$ 
such that $g(\xi_i)=sign(\varphi_i(x^\ast))$, $1\leq i\leq M$, 
since $g$ is trivially continuous on the closed discrete subset $\Xi\subset K$, there exists an extension $g^\ast\in C(K)$ taking values in $[-1,1]$, with $L_n g^\ast(x^\ast)=\lambda_n(x^\ast)$.

Now, applying Lemma 1 to the polynomial $L_nf$ we get 
$$
\|L_nf\|_K\leq c_m \|L_nf\|_{X_n^m}\;.
$$
On the other hand $|L_nf(x)|\leq \|f\|_\Xi\lambda_n(x)\leq \|f\|_K \lambda_n(x)$ and thus $\|L_nf\|_K\leq c_m \|f\|_K \|\lambda_n\|_{X_n^m}$, from which we get immediately 
$$
\|L_n\|=\|\lambda_n\|_K 
\leq c_m\|\lambda_n\|_{X_n^m}\;,
$$
and thus (\ref{lebest}), since $\|\lambda_n\|_K\geq \|\lambda_n\|_{X_n^m}$  by inclusion.
\hspace{0.2cm} $\square$
\vskip0.3cm

\begin{Remark}
{\em 
Notice that $c_m\to 1$ and thus, if the sampling set $\Xi$ is independent of $m$, $\|\lambda_n\|_{X_n^m}\to \|L_n\|$ as $m\to \infty$. Moreover, from (\ref{lebest}) we also get a relative error estimate, 
namely
\begin{equation} \label{asymptotic}
0\leq \frac{\|L_n\|
-\|\lambda_n\|_{X_n^m}}{\|L_n\|}
\leq (c_m-1)\sim \frac{\pi^2}{8m^2}\approx \frac{1.23}{m^2}\;, 
\end{equation}
i.e. $\|\lambda_n\|_{X_n^m}$ approximates the Lebesgue constant 
from below with a 
$\mathcal{O}(1/m^2)$ relative error.
On the other hand, (\ref{lebest}) gives also 
the rigorous and {\em computable\/} absolute error estimate $0\leq \|L_n\| 
-\|\lambda_n\|_{X_n^m}\leq (c_m-1) \|\lambda_n\|_{X_n^m}$.
We finally notice that (\ref{lebest}) is an {\em interval approximation} of the Lebesgue constant. Hence we can use the midpoint approximation 
\begin{equation} \label{midpoint}
\frac{|\,\|L\|_n-\Lambda_n^m\,|}{\|L_n\|}\leq \frac{c_m-1}{2}\;,\;\;\Lambda_n^m=\|\lambda_n\|_{X_n^m}\,\frac{1+c_m}{2}\;,
\end{equation}
which improves the error estimates by a factor $1/2$.  
\/}
\end{Remark}

\begin{Remark}
{\em 
The structure of projection operators like (\ref{Ln}) includes interpolation operators at unisolvent nodes $\Xi=\{\xi_1,\dots,\xi_{N}\}\subset K$, where, denoting by 
$V_n=[p_j(\xi_i)]$, $1\leq i,j \leq N$, the Vandermonde-like matrix in any fixed polynomial basis $span\{p_1,\dots,p_N\}=\mathbb{P}_n$, we have that 
\begin{equation} \label{interp}
\varphi_j(x)=\ell_j(x)=\frac{det(V_n(\xi_1,\dots,\xi_{j-1},x,\xi_{j+1},\dots,\xi_N))}{det(V_n(\xi_1,\dots,\xi_{j-1},\xi_j,\xi_{j+1},\dots,\xi_N))}
\end{equation}
are the corresponding Lagrange cardinal polynomials. But also discrete weighted least-squares operators at  $\mathbb{P}_n$-determining nodes $\Xi=\{\xi_1,\dots,\xi_M\}\subset K$ with positive weights 
$W=\{w_1,\dots,w_M\}$, $M>N$, are included. Indeed, denoting by $\{\pi_k\}$, $1\leq k\leq N$, the orthonormal polynomials with respect to the corresponding discrete scalar product $(f,g)_{\ell^2 _W(\Xi)}=\sum_{j=1}^M{w_jf(\xi_j)g(\xi_j)}$, we have that
\begin{equation} \label{ls}
L_nf(x)=\sum_{k=1}^{N}{(f,\pi_k)_{\ell^2 _W(\Xi)}\,\pi_k(x)}=
\sum_{j=1}^M{f(\xi_j)w_jK_n(x,\xi_j)}\;,
\end{equation}
i.e. $\varphi_j(x)=w_jK_n(x,\xi_j)$, where $K_n(x,y)=\sum_{k=1}^{N}{\pi_k(x)\pi_k(y)}$ is the reproducing kernel of the discrete scalar product. Notice that in this case (unless $M=N$ where least-squares approximation coincides with interpolation) the $\varphi_j$ are linearly dependent, thus forming a set of generators of $\mathbb{P}_n$.

We stress that the structure (\ref{ls}) also includes {\em hyperinterpolation} operators, a topic that has seen an increasing interest as a valid alternative to multivariate polynomial interpolation, after the seminal paper \cite{Sl95} by Sloan in the mid '90s; cf., e.g., \cite{AW21,DMVX09,HAC09,SH07,SV17,Wade13,WWW14} with the references therein. Indeed hyperinterpolation operators are substantially truncated Fourier-like expansions in series of orthogonal polynomials with respect to a continuous measure with density. The continuous scalar products are there substituted by discrete scalar products, corresponding to a suitable positive quadrature formula exact in $\mathbb{P}_{2n}$. 
\/}
\end{Remark}

\subsection{The simplex}

We consider now the $d$-simplex
\begin{equation} \label{simplex}
T_d=\{x=(x_1,\dots,x_d)\in \mathbb{R}^d:\,0\leq x_d\leq \dots\leq x_1\leq 1\}\;,
\end{equation}
along with the $d$-dimensional Duffy-like transformation $\mathcal{D}:\,[-1,1]^d\to T_d$, 
which can be defined as follows (cf. e.g. \cite{ST23})
\begin{equation} \label{duffy}
x_i=\mathcal{D}_i(t)=\prod_{j=1}^i(t_j/2+1/2)\;,\;
1\leq i\leq d\;\;,\;\;t=(t_1,\dots,t_d)\in [-1,1]^d\;.
\end{equation}
We again notice that the results can be immediately extended to any simplex by invertible affine transformation, 
with invariance of the mesh constant $c$.

\begin{Lemma}
Let $\mathcal{C}_k$ be the Chebyshev or Chebyshev-Lobatto points in $[-1,1]$ and $\mathcal{D}$ the Duffy-like transformation $[-1,1]^d\to T_d$ in (\ref{duffy}), where $T_d$ is the $d$-simplex. Then the sequence $X_n^m=\mathcal{D}((\mathcal{C}_{mn})^d)$, $n=1,2,\dots$, for a  fixed $m>1$, is an admissible polynomial mesh for $T_d$, with constant $c=(c_m)^d$.
\end{Lemma}
\vskip0.3cm
\noindent{\em Proof.} 
Let us denote by $\mathbb{P}_n^1$ the space of univariate real algebraic polynomials of degree not exceeding $n$. For every $p\in \mathbb{P}_n$ we have that $\|p\|_{T_d}=\|p\circ \mathcal{D}\|_{[-1,1]^d}$. Now, since $\mathcal{D}$ is a $d$-linear (surjective) mapping, $p\circ \mathcal{D}\in \bigotimes_{k=1}^d\mathbb{P}_n^1$, the space of tensorial polynomials of degree not exceeding $n$. Reasoning component by component and using iteratively the univariate version of Lemma 1, it is immediate to write the inequality $$\|p\circ \mathcal{D}\|_{[-1,1]^d}\leq (c_m)^d\|p\circ \mathcal{D}\|_{(\mathcal{C}_{mn})^d}=(c_m)^d\|p\|_{X_n^m}\;.\hspace{1cm}\square
$$

\begin{Proposition}
Let $X_n^m=\mathcal{D}((\mathcal{C}_{mn})^d)$ be a polynomial admissible mesh of the $d$-simplex $K=T_d$ as in Lemma 2, and $L_n:C(K)\to \mathbb{P}_n$ a linear projection operator with the structure defined in Proposition 1. 
Then, the following estimate holds for every $m>1$ 
\begin{equation} \label{lebest2}
\|\lambda_n\|_{X_n^m}\leq \|L_n\|\leq 
(c_m)^d\|\lambda_n\|_{X_n^m}\;.
\end{equation} 

\end{Proposition}
\vskip0.3cm
\noindent{\em Proof.} We can proceed exactly as in the proof of Proposition 1, simply by substituting $c_m$ with $(c_m)^d$.
\hspace{0.2cm} $\square$

\begin{Remark}
{\em 
We observe that, since $(c_m)^d\to 1$ as $m \to \infty$, the same assertions of Remark 1 are valid, with $(c_m)^d$ replacing $c_m$. The only difference is that estimate (\ref{asymptotic}) is asymptotically increased by a factor $d$, namely 
$$
0\leq \frac{\|L_n\|
-\|\lambda_n\|_{X_n^m}}{\|L_n\|}\leq (c_m)^d-1
=(c_m^{d-1}+\dots+c_m+1)(c_m-1)
$$
\begin{equation} \label{asymptotic2}
\leq \left(\left(\frac{m}{m-1}\right)^{d-1}+\dots+\frac{m}{m-1}+1\right)(c_m-1)
\sim \frac{d\pi^2}{8m^2}\approx 1.23\,\frac{d}{m^2}
\end{equation}
for $d$ fixed and $m\to \infty$, where we have used the 
elementary inequality $\cos(\theta)=\sin(\pi/2-\theta)\geq 1-2\theta/\pi$, $0\leq \theta \leq \pi/2$.
}
\end{Remark} 

\begin{Remark}
{\em 
It is worth observing that, by the finite union extension property of admissible polynomial meshes (cf. e.g. \cite{CL08}), the results above are valid via ``triangulation'' (subdivision into non overlapping simplices) on any polytope, e.g. on any polygon in $d=2$ or polyhedron in $d=3$. We recall that for $d\geq 3$ in nonconvex instances the triangulation could require extra vertices, in view of the well-known Sch\"{o}nardt counterexample \cite{S28}. 
On the other hand, for the same reason it is also valid on any union of (even overlapping) polytopes, where again the mesh is the union of the single polytope meshes obtained via triangulation (in practice, this avoids tracking and triangulating the union polytope which can be very complicated). 
}
\end{Remark} 

\subsection{The ball}
As a third relevant case, we discuss the unit Euclidean ball, i.e.  $B_d=\{x\in \mathbb{R}^d:\,\|x\|_2\leq 1\}$, by no loss of generality since it is affinely equivalent to any other ball by translation and scaling, with invariance of the mesh constant $c$. Again, the results can be extended to a finite union of possibly overlapping balls, a geometrical object that is relevant in applications, e.g. 
in the field of molecular modelling \cite{P13}.

Below we denote by $\mathbb{P}_n^1$ the space of univariate real algebraic polynomials of degree not exceeding $n$, and by $\mathbb{T}_n^1([a,b])$ the space of univariate real trigonometric 
polynomials of degree not exceeding $n$, i.e $span\{1,\sin(j\theta),\cos(j\theta)\,,\,j=1,\dots,n\}$, restricted to a subinterval $[a,b]$ of the period, $b-a\leq 2\pi$. 

Moreover, we make use of the {\em generalized spherical coordinates} in $B_d$, $d\geq 2$ ($B_1=[-1,1]$ is  treated in section 2.1), corresponding to the surjective transformation $$\mathcal{G}: J=[0,1]\times [0,\pi]^{d-2}\times [0,2\pi]\to B_d$$ defined by 
$$
\;x_j=r\cos(\theta_j)\prod_{k=1}^{j-1}\sin(\theta_k)\;,\;1\leq j\leq d-2\;,
$$
\begin{equation} \label{sphcoord}
x_{d-1}=r\cos(\theta_{d-1})\prod_{k=1}^{d-2}\sin(\theta_k)\;,
\;x_d=r\sin(\theta_{d-1})\prod_{k=1}^{d-2}\sin(\theta_k)\;,
\end{equation}
where$\;r\in [0,1]\;,\;\;\theta_k\in [0,\pi]\;,\;1\leq k\leq d-2$, 
$\theta_{d-1}\in [0,2\pi]$; cf., e.g., \cite{Blu60}. These coordinates coincide with the usual polar coordinates for the disk $B_2$ and spherical coordinates for the 3-ball $B_3$.

First, we state a basic Lemma on a norming inequality for 
univariate trigonometric polynomials in the subperiodic case, whose proof can be found in \cite{V18-2}.

\begin{Lemma}
Let $\mathcal{C}_k$ be the Chebyshev or Chebyshev-Lobatto points in $[-1,1]$ and let $\sigma_{a,b}:[-1,1]\to[a,b]$, $b-a\leq 2\pi$, be the invertible map
\begin{equation} \label{sigmatrig}
\sigma_{a,b}(u)=2\arcsin(\alpha u)+\beta\;,\;\;u\in [-1,1]\;,\;
\alpha=\sin\left(\frac{b-a}{4}\right)\;,\;\beta=\frac{b+a}{2}\;.
\end{equation}
Then the following inequality holds
\begin{equation} \label{trigmesh}
\|\phi\|_{a,b}\leq c_m \|\phi\|_{\sigma_{a,b}(\mathcal{C}_{2mn})}\;,\;\forall \phi\in \mathbb{T}_n([a,b])\;.
\end{equation}
\end{Lemma}

\begin{Remark}
{\em We observe that the Chebyshev-like angles $\sigma_{a,b}(\mathcal{C}_{2mn})$ cluster at the interval endpoints for $b-a<2\pi$ (subperiodic instances), whereas are equally spaced for $b-a=2\pi$ (periodic case).
}
\end{Remark}

\begin{Lemma}
Let $\mathcal{C}_k$ be the Chebyshev or Chebyshev-Lobatto points in $[-1,1]$ and consider the composed transformation 
\begin{equation} \label{comptrasf}
\mathcal{S}=\mathcal{G}\circ \mathcal{U}:[-1,1]^d\to B_d\;,
\end{equation}
$$
\mathcal{U}(u_1,u_2,\dots,u_{d-1},u_d)=\left(\frac{u_1}{2}+\frac{1}{2},\sigma_{0,\pi}(u_2),\dots,\sigma_{0,\pi}(u_{d-1}),\sigma_{0,2\pi}(u_d)\right)\;,
$$
where $B_d$ is the $d$-ball. Then the sequence $X_n^m=\mathcal{S}(\mathcal{C}_{mn}\times (\mathcal{C}_{2mn})^{d-1})$, $n=1,2,\dots$, for a  fixed $m>1$, is an admissible polynomial mesh for $B_d$, with constant $c=(c_m)^d$.

\end{Lemma}
\vskip0.3cm
\noindent{\em Proof.} For every $p\in \mathbb{P}_n$, it is easily seen by basic trigonometric identities that the composed function $p\circ \mathcal{G}$ belongs to an algebraic-trigonometric tensorial space 
on the box $J=[0,1]\times [0,\pi]^{d-1}\times [0,2\pi]=\mathcal{U}([-1,1]^ d)$, namely 
$$p\circ \mathcal{G}\in \mathbb{P}_n^1\bigotimes \mathbb{T}_n([0,\pi])\bigotimes \dots \bigotimes \mathbb{T}_n([0,\pi]) \bigotimes \mathbb{T}_n([0,2\pi])\;.$$ 
Moreover, $\|p\|_{B_d}=\|p\circ \mathcal{G}\|_{J}$ by surjectivity of $\mathcal{G}$. 
Then, reasoning component by component by using the univariate version of Lemma 1 and iteratively Lemma 3, it is immediate to write $$\|p\|_{B_d}=\|p\circ \mathcal{G}\|_{J}\leq (c_m)^d\|p\circ \mathcal{G}\|_{\mathcal{U}(\mathcal{C}_{mn}\times (\mathcal{C}_{2mn})^{d-1})}
$$
$$
=(c_m)^d\|p\circ \mathcal{S}\|_{\mathcal{C}_{mn}\times (\mathcal{C}_{2mn})^{d-1}}=(c_m)^d\|p\|_{X_n^m}\;.
\hspace{1cm} \square
$$

\begin{Proposition}
Let $X_n^m=\mathcal{S}(\mathcal{C}_{mn}\times (\mathcal{C}_{2mn})^{d-1})$ be a polynomial admissible mesh of the $d$-ball $K=B_d$ as in Lemma 4, and $L_n:C(K)\to \mathbb{P}_n$ a linear projection operator with the structure defined in Proposition 1. 
Then, the following estimate holds for every $m>1$ 
\begin{equation} \label{lebest2}
\|\lambda_n\|_{X_n^m}\leq \|L_n\|\leq 
(c_m)^d\|\lambda_n\|_{X_n^m}\;.
\end{equation}

\end{Proposition}
\vskip0.3cm
\noindent{\em Proof.} Again (see Proposition 2), we can proceed exactly as in the proof of Proposition 1, simply by substituting $c_m$ with $(c_m)^d$.\hspace{0.2cm} $\square$
\vskip0.3cm

We can finally observe that in view of (\ref{lebest2}), the considerations in Remark 3 apply also to the case of the ball.

\subsection{General polynomial meshes}

In the case of more general compact sets, we can still get a fully discrete but less accurate approximation of Lebesgue constants by polynomial meshes, based on the approach developed in \cite{V18} for polynomial optimization. 

\begin{Proposition}
Let $K\subset \mathbb{R}^d$ be a compact set, $\{X_n\}$ a polynomial admissible mesh of $K$, and $L_n:C(K)\to \mathbb{P}_n$ a linear projection operator with the structure defined in Proposition 1. 
Then, the following estimates hold for every $m\geq 1$ 
\begin{equation} \label{lebest2}
\|\lambda_n\|_{X_{mn}}\leq \|L_n\|\leq 
c^{1/m}\|\lambda_n\|_{X_{mn}}\;.
\end{equation} 

\end{Proposition}
\vskip0.3cm
\noindent{\em Proof.} First, observe that by definition of polynomial mesh, for every $p\in \mathbb{P}_n$ and for every $m\geq 1$ we can write
$$
\|p\|_K\leq c^{1/m} \|p\|_{X_{mn}}\;,
$$
since $p^m\in \mathbb{P}_{mn}$ and $\|p^m\|_K=\|p\|_K^m\leq C \|p^m\|_{X_{mn}}
=C\|p\|_{X_{mn}}^m$. 
Then we can reason as in the proof of Proposition 1 to reach the conclusion, 
with $c^{1/m}$ substituting $c_m$. 
\hspace{0.2cm} $\square$
\vskip0.3cm

\begin{Remark}
{\em  
Notice that $c^{1/m}\to 1$ and thus again, if the sampling set $\Xi$ is independent of $m$, $\|\lambda_n\|_{X_{mn}}\to \|L_n\|$ as $m\to \infty$. Moreover, from (\ref{lebest2}) we also get a relative error estimate, 
namely
\begin{equation} \label{asymptotic2}
\frac{\|L_n\| -\|\lambda_n\|_{X_{mn}}}{\|L_n\|} 
\leq e^{\log(c)/m}-1\leq e^{\log(c)/m}\, \frac{\log(c)}{m}
\end{equation}
by
the mean value theorem and the monotonicity of the exponential function, that is a $\mathcal{O}(1/m)$ relative approximation of the Lebesgue constant by $\|\lambda_n\|_{X_{mn}}$. Again,  (\ref{lebest2}) gives also 
the rigorous and {\em computable\/} absolute error estimate $0\leq \|L_n\| 
-\|\lambda_n\|_{X_{mn}}\leq (c^{1/m}-1) \|\lambda_n\|_{X_{mn}}$. Notice that also in this general case we may resort to the midpoint approximation 
$\Lambda_{X_{mn}}=\|\lambda_n\|_{X_{mn}}(1+c^{1/m})/2$ whose relative error estimate is improved by a factor $1/2$, namely 
$|\,\|L_n\| 
-\Lambda_{X_{mn}}\,|/\|L_n\|\leq (c^{1/m}-1)/2$}.
\end{Remark}

\section{Numerical examples}

\subsection{Computational issues} 
We make some observations on the main computational issues. Dealing with polynomial projectors, 
most computations can be seen as a matter of numerical linear algebra, by standard algorithms applied to the relevant Vandermonde-like matrices. 
Indeed, given a finite set of points $Z=\{z_1,\dots,z_k\}$, and selected a polynomial basis $[p_1(x),\dots,p_{N}(x)]$ of $\mathbb{P}_n$, we can construct the matrix 
\begin{equation} \label{vand}
V_n(Z)=[p_j(z_i)]\;,\;1\leq i\leq k\;,\;1\leq j\leq N\;.
\end{equation}
In order to control the conditioning of such matrices, that becomes unacceptable already at moderate degrees with the standard monomial basis, we have chosen to use a total-degree product Chebyshev basis corresponding to the minimal enclosing box, say 
$[a_1,b_1]\times\dots\times[a_d,b_d]$, namely $p_j(x)=\prod_{s=1}^dT_{m_s}(\alpha_sx_s+\beta_s)$ where $\alpha_s=(b_s-a_s)/2$, $\beta_s=(b_s+a_s)/2$ and $T_{m_s}(\cdot)=\cos(m_s\arccos(\cdot))$ is the standard Chebyshev polynomial of the second-kind for degree $m_s$. Here $j=1,\dots,N$ corresponds to some ordering (for example a lexicographical ordering) of the $d$-uples $(m_1,\dots,m_d)$, $0\leq m_1+...+m_d\leq n$. 

Alternatively, when known one could adopt orthonormal polynomial bases with respect to some continuous measure with density  on the specific domain, such as the Logan-Shepp or the Zernike basis for the disk \cite{W08}, the Dubiner basis for the 2-simplex \cite{Dub91}.

Notice that, in view of (\ref{ls}), the Lebesgue constant on the mesh of either an interpolation or weighted least-squares projection operator, can then be simply computed via the relevant matrices as 
\begin{equation} \label{matrixleb}
\|\lambda_n\|_{X_n^m}=\max_i\sum_j w_j\left|\sum_k \pi_k(x_i)\pi_k(\xi_j)\right|
=\|V_n(X_n^m)R^{-1}Q^tdiag(\sqrt{W})\|_\infty\;,
\end{equation}
 where $\Xi=\{\xi_j\}$ are either the interpolation or the weighted least-squares sampling points, and 
\begin{equation} \label{qr}
diag(\sqrt{W})\,V_n(\Xi)=QR
\end{equation}
the factorization of the corresponding weighted Vandermonde-like matrix, with $Q$ (rectangular) orthogonal and $R$ square upper-triangular, that is 
\begin{equation} \label{orth}
[\pi_1(x),\dots,\pi_{N}(x)]=[p_1(x),\dots,p_{N}(x)]R^{-1}
\end{equation}
is the corresponding discrete orthonormal polynomial basis in $\ell^2_W(\Xi)$.
\vskip0.5cm

To have an idea of the approximation quality, we report here a table of the relative error estimates using the midpoint approximation (cf. (\ref{midpoint}), with $c_m^d$ replacing $c_m$ for simplex and ball). Observe that to get an error below 10\%, that is to recover the Lebesgue constant with one correct figure thus computing accurately its order of magnitude (which is the relevant parameter in applications), one can take $m=3$ for the $d$-cube, $m=4$ in $d=2$ for simplex and disk, and $m=5$ in $d=3$ for simplex and ball (in all cases the error 
is around $8\%$).

\begin{table}[ht]
\begin{center}
{\footnotesize
\begin{tabular}{|c|c|c|c|c|c|c|c|c|c|}
\hline
$m$ & 2 & 3 & 4 & 5 & 6 & 7 & 8 & 9 & 10\\
\hline \hline
$c_m-1$ & 11\% & 7.7\% & 4.1\% & 2.6\% & 1.8\% & 1.3\% & 0.98\% & 0.77 \% & 0.62\% \\ 
\hline
$c_m^2-1$ & 50\% & 17\% & 8.6\% & 5.3\% & 3.6\% & 2.6\% & 2.0\% & 1.6\% & 1.3\% \\ 
\hline
$c_m^3-1$ & 90\% & 27\% & 13\% & 8.1\% & 5.5\% & 4.0\% & 3.0\% & 2.4\% & 1.9\% \\ 
\hline
\end{tabular}
}
\caption{\footnotesize{Relative error estimates in the approximation of the Lebesgue constants on Chebyshev polynomial meshes $X_n^m$: 
$d$-cube (first row), 2-simplex and disk (second row), 3-simplex and ball (third row).}}
\label{relerrors}
\end{center}
\end{table}

\begin{figure}[h!]
    \centering
        {\includegraphics[scale=0.4]{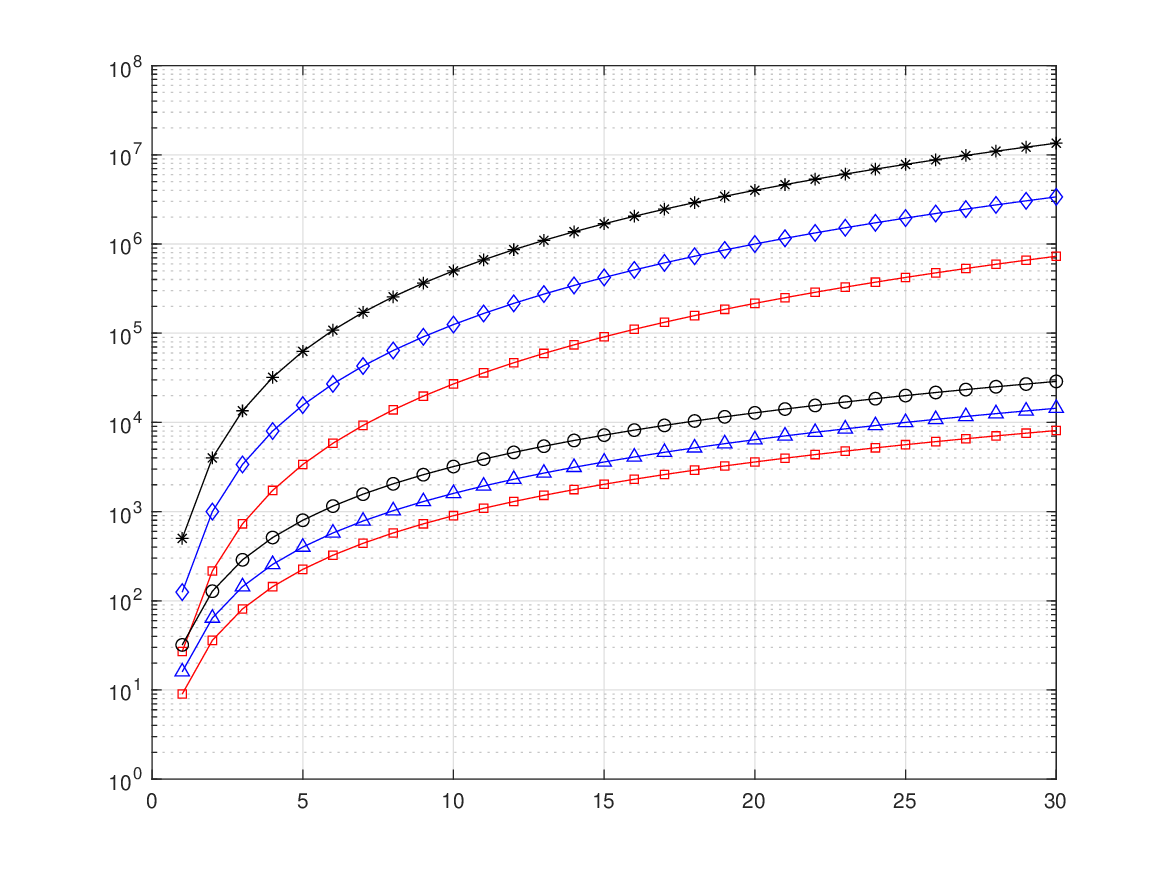}}
    \caption{\footnotesize{
    Cardinalities (log scale) of some Chebyshev polynomial meshes in dimension $d=2,3$ for degree $n=1,2,\dots,30$ (Lebesgue constant approximated at less than $10\%$, cf. Table 1): $(mn)^d$ with $m=3$ ($d$-cube, $\square$), $m=4$ (2-simplex, $\triangle$) and $m=5$ (3-simplex, $\lozenge$); $2^{d-1}(mn)^d$ with $m=4$ (disk, o) and $m=5$ (ball, $\ast$).}}
    \label{card}
\end{figure}

We observe that the present fully discrete approach for Lebesgue constants evaluation, being based on product Chebyshev meshes of cardinality $\mathcal{O}((mn)^d)$, suffers of the curse of dimensionality and hence is essentially a low-dimensional tool. 
Considering for example the polynomial meshes corresponding to (transformed) grids of Chebyshev points, which are in the interior 
of the domains, the cardinalities are exactly $(mn)^d$ for the $d$-cube and the $d$-simplex, and  $2^{d-1}(mn)^d$ for the $d$-ball.
To have an idea of the sizes, we plot in Figure \ref{card} the values of the cardinality corresponding to the choices of $m$ suggested above, for a range of polynomial degrees in dimension $d=2$ and $d=3$.

\vskip0.2cm

\vskip0.3cm 

Below we show a number of numerical tests in dimension 1, 2 and 3. The corresponding numerical codes and demos, implemented in Matlab, are freely available at \cite{KSV23}.

\subsection{Univariate interpolation points}
\begin{figure}[h!]
    \centering
        {\includegraphics[scale=0.3]{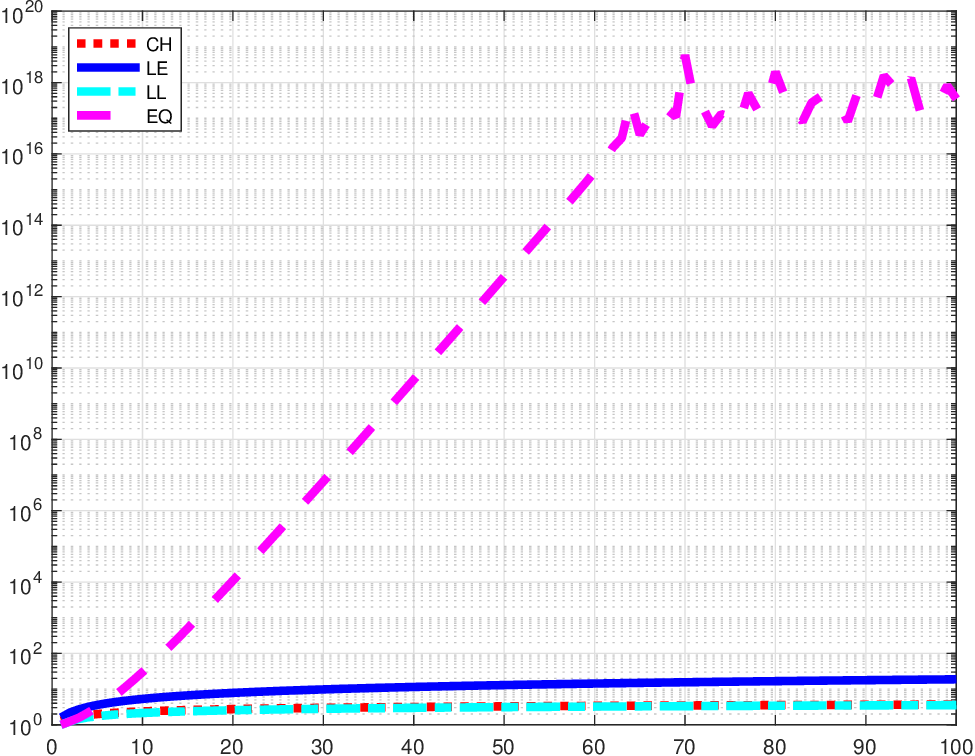}}{\includegraphics[scale=0.3]{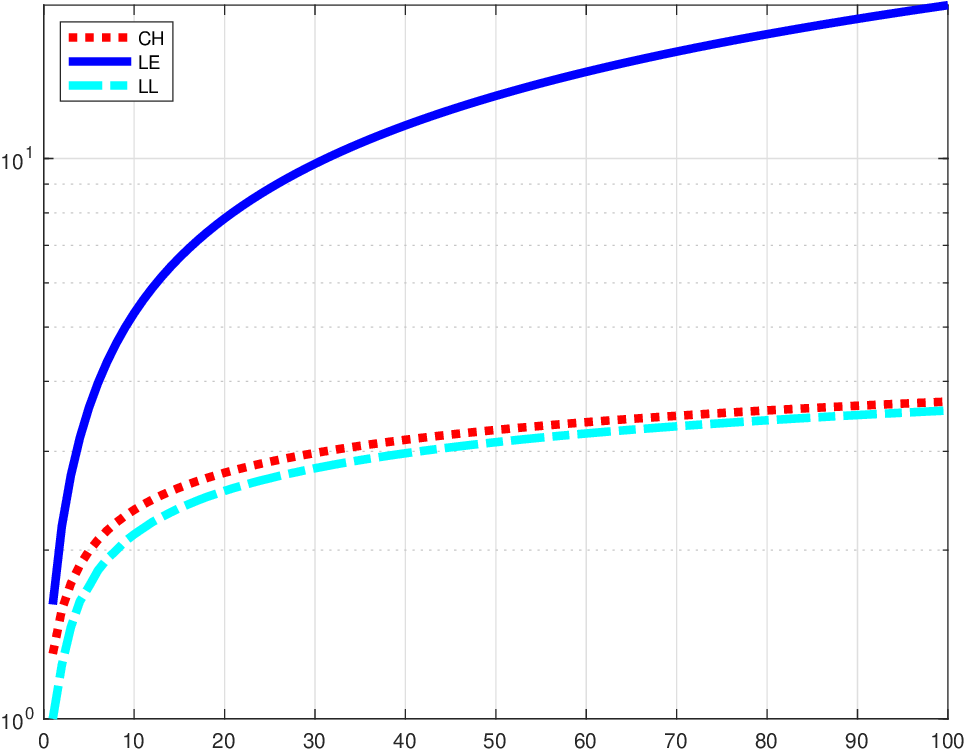}}
    \caption{\footnotesize{
    Left: Lebesgue constants of some pointsets on the interval $[-1,1]$  for degrees $n=1,\dots,100$: Chebyshev (red dots), Legendre (blue line), Legendre-Lobatto (cyan dashes), equispaced points (magenta dashes). 
    Right: the three lowest curves in detail.
    In these experiments, $m=3$ with a relative error $\approx 7.7\%$. }}
    \label{leb_intv}
\end{figure}

In the univariate case Lebesgue constants for interpolation on intervals have been extensively studied, with a number of theoretical results and estimates; cf. e.g. \cite{Brut97,MM08} and the references therein. In order to test our method, we compare here 
the computed Lebesgue constants of Chebyshev points, Gauss-Legendre-Lobatto points (which are known to be Fekete points, i.e. points that maximize the absolute value of the Vandermonde determinant) and Gauss-Legendre points. The Lebesgue constant of the first two is known to be $\mathcal{O}(\log(n))$, whereas the third is $\mathcal{O}(\sqrt{n})$. Moreover, we also compute the Lebesgue constant of equally spaced points, which is known to grow exponentially (cf. \cite{MS92}). The results are collected in Figure \ref{leb_intv}.  

\vskip0.2cm

\subsection{Comparing interpolation points on the square}
\begin{figure}[h!]
    \centering
        {\includegraphics[scale=0.3]{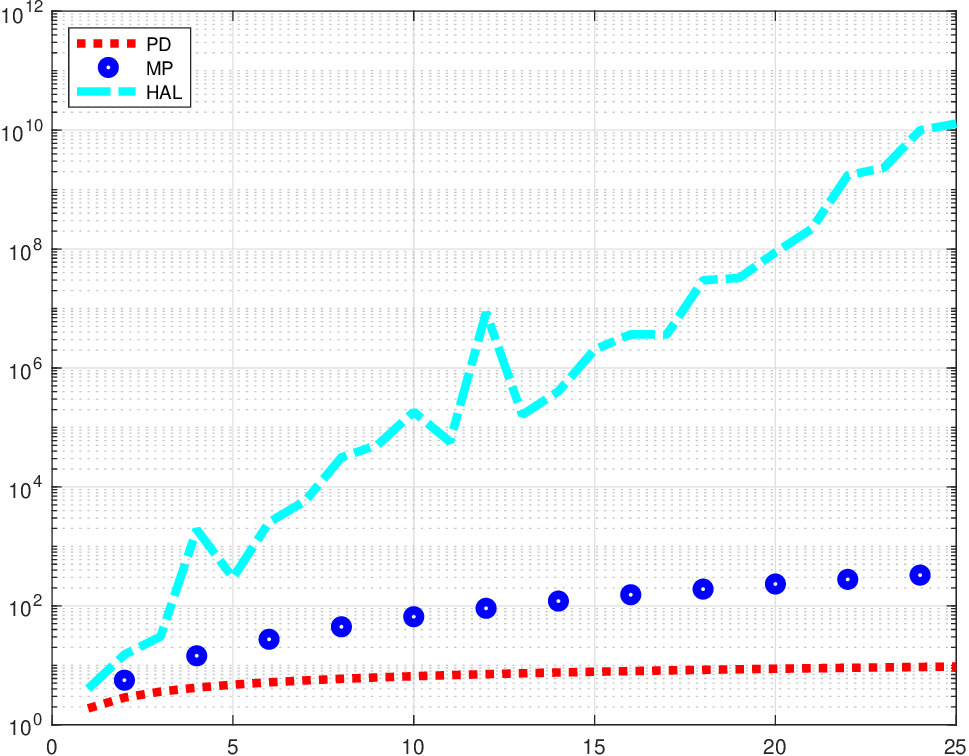}}
        {\includegraphics[scale=0.3]{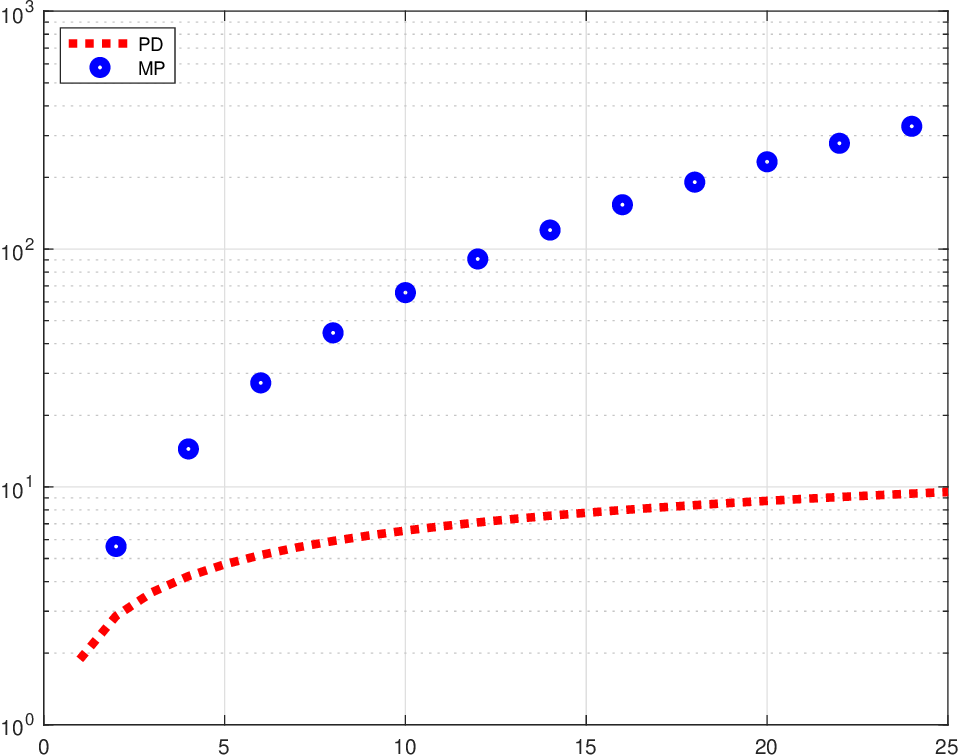}}
    \caption{\footnotesize{
    Left: Lebesgue constants of some pointsets on the square $[-1,1]^2$  for degrees $n=1,\dots,25$: Padua points (red dots), Morrow-Patterson (blue dots), Halton points (cyan dashes). 
    Right: the two lowest curves in detail.
    In these experiments, $m=3$ with a relative error $\approx 7.7\%$.}}
    \label{leb_sq}
\end{figure}

We compare here the computed Lebesgue constants of some well-known families of interpolation points on the square. The results are collected in Figure \ref{leb_sq}. 

The Padua points on the square, 
discovered in 2005 \cite{CDMV05}, are the union of two suitable Chebyshev subgrids. They are the first and till now the only explicitly known optimal point set for total-degree multivariate polynomial interpolation. For such points it has indeed been proved that the Lebesgue
constant is $\mathcal{O}(\log^ 2(n))$; cf.  \cite{BCDMVX06}. 

The Morrow-Patterson points support one of the few known minimal positive cubature formulas, namely a
formula with $N$ nodes that has degree of exactness $2n$ for the product Chebyshev measure of the second kind, cf. \cite{MP78}. Hence, the hyperinterpolation polynomial of degree not greater than $n$ at these points, in view of
minimality, turns out to be exactly the interpolation polynomial, by \cite[Lemma 3]{Sl95}. For such points it is proved that the Lebesgue constant is $\mathcal{O}(n^3)$ and it is conjectured that the actual growth is $\mathcal{O}(n^2)$; cf. \cite{DMSV14}.

For the  purpose of comparison we also compute the Lebesgue constant of $N$ Halton points. In view of the recent result \cite{DASV23}, $N$ uniformly distributed (random) points are almost surely unisolvent, the probability that $det(V_n(\Xi))=0$ being null, but we expect that the Lebesgue constant has exponential growth, as observed numerically. A theoretical explanation is that a subexponential growth would imply weak-* convergence of the uniform discrete probability measure supported at the interpolation points, to the potential theoretic equilibrium measure of the compact set (that in this case is the product Chebyshev measure); cf. \cite{BBCL12}. On the contrary, with uniform random and Halton points there is weak-* convergence to the Lebesgue measure (a fact that is at the base of MonteCarlo and Quasi-MonteCarlo integration). The exponential growth is clearly visible in Figure \ref{leb_sq}-left, and in all the figures with uniform or Halton points as a trend up to some oscillations (in Figure \ref{leb_intv}-left where high degrees are considered, the behavior becomes numerically erratic when the Lebesgue constant goes beyond a size  around $10^{17}$). 

\vskip0.2cm

\subsection{Waldron points on the simplex}

Good total-degree interpolation points 
on the simplex are relevant in the numerical solution of PDEs by spectral and high-order methods, and for this reason have been extensively investigated, in most cases numerically; cf. e.g. \cite{W06} with the references therein. 

\begin{figure}[h!]
    \centering
        {\includegraphics[scale=0.3]{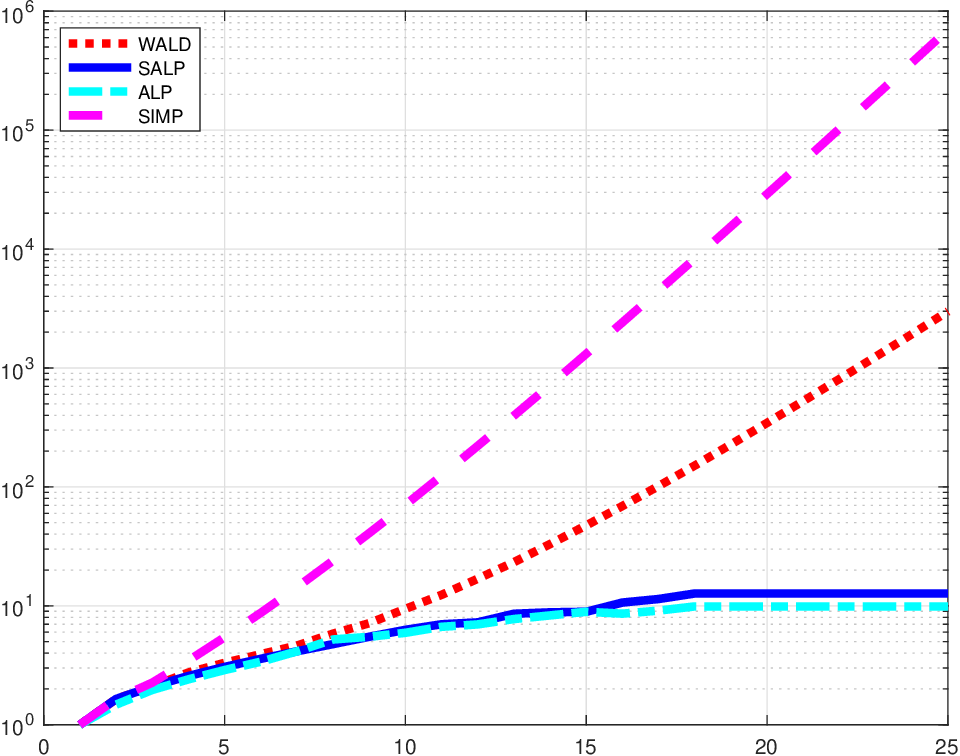}}
        {\includegraphics[scale=0.3]{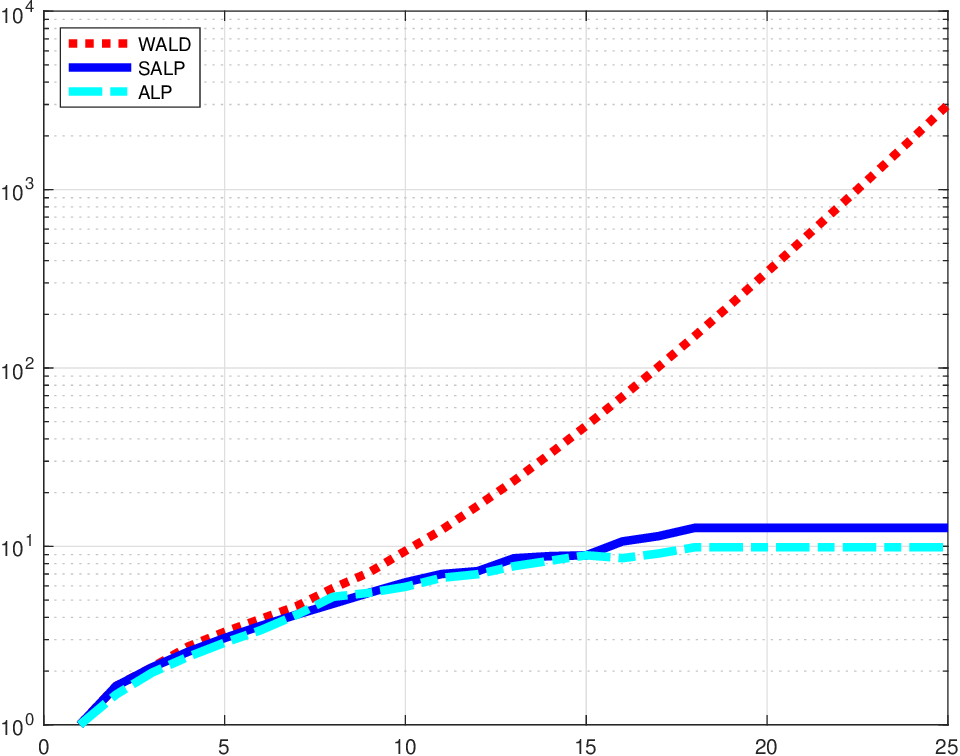}}
    \caption{\footnotesize{
    Left: Lebesgue constants of some pointsets on the unit simplex with vertices $(0,0)$, $(1,0)$, $(0,1)$, for degrees $n=1,\dots,25$: ALP (cyan dashes), SALP (blue line), Waldron points (red dots), simplex points (magenta dashes). 
    Right: the three lowest curves in detail. 
    In these experiments, $m=4$ with a relative error $\approx 8.6\%$.}}
    \label{leb_sym}
\end{figure}

Quite recently a promising theoretical approach has been proposed, constructing the so-called Waldron points, that are obtained by looking for an appropriate spacing with respect to a distance related to the equilibrium measure of the domain (such as the Baran distance; cf. \cite{BMW23}). In Figure \ref{leb_sym} we compare the Lebesgue constant of the Waldron points for the 2-simplex with that of the so-called Simplex Points (SIMP) (a triangular grid corresponding to equally spaced points in the Euclidean distance on the equilateral triangle), and two families of points corresponding to a greedy minimization of the Lebesgue constant, the Approximate Lebesgue Points (ALP) and the Symmetric Approximate Lebesgue Points (SALP). The SALP are useful in the framework of spectral element methods for PDEs; cf. \cite{BSV12,RSV12}.

We see that as expected with the Simplex Points there is an exponential growth. On the other hand, the Lebesgue constant of ALP and SALP are slowly increasing, whereas that of the Waldron points increase slowly up to about degree $n=10$ and then it turns to a manifestly exponential growth (though slower than with the Simplex Points).

\vskip0.2cm

\subsection{Approximate Fekete and Leja points}
The cases where good interpolation sets are known analytically  
are very few. Already in dimension $d=2$ to our knowledge there is no known explicit family for the disk, and in $d=3$ the same can be said for cube and ball. 

On the other hand, interpolation points with slowly increasing Lebesgue constant can be determined numerically, for example the Approximate Fekete Points (AFP) and the Discrete Leja Points (DLP), both corresponding to a greedy maximization of the Vandermonde determinant modulus, typically extracting such points from polynomial meshes by numerical linear algebra algorithms; cf., e.g., \cite{BCLSV11,BDMSV10}. 
Moreover, for example in \cite{BSV12,MS19} Approximate Lebesgue Points (ALP) have been computed once and forall on square, simplex and disk for restricted degree ranges, working heuristically just with polynomial meshes and suitable greedy algorithms. Other pointsets with low Lebesgue constant can  be computed by the optimization algorithm proposed in \cite{VBHS14}.

\begin{figure}[h!]
    \centering
        {\includegraphics[scale=0.3]{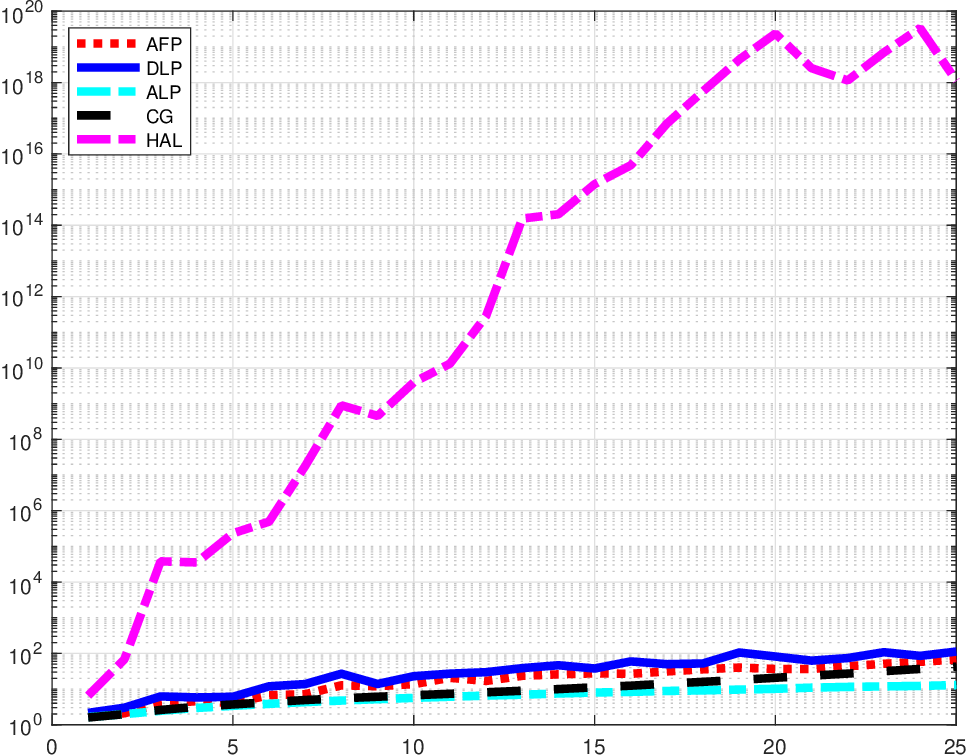}}
        {\includegraphics[scale=0.3]{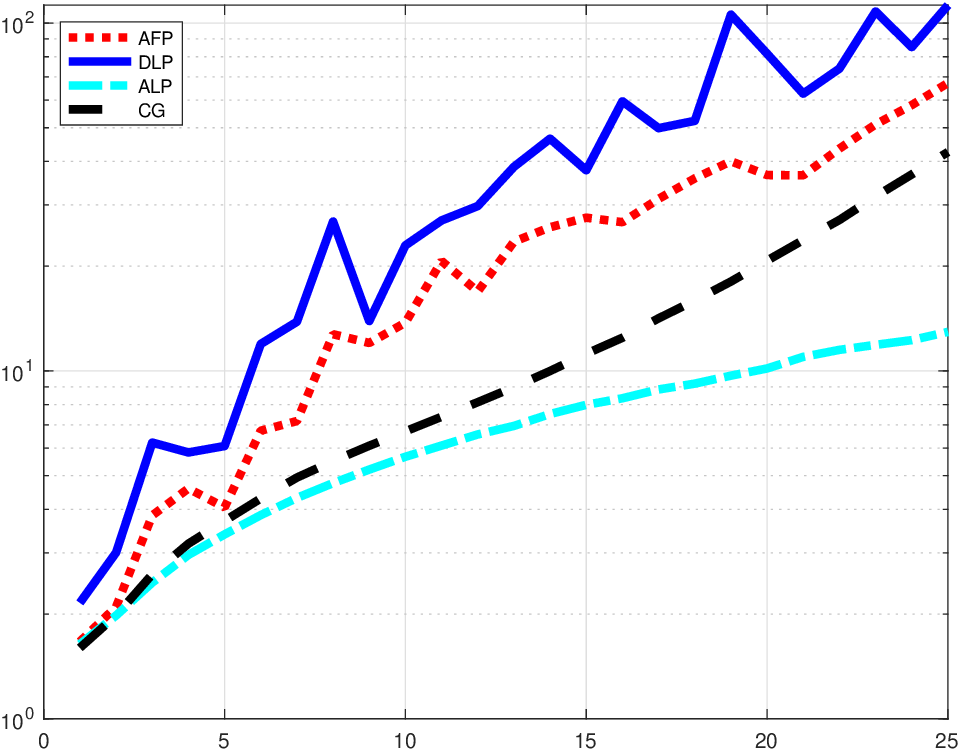}}
    \caption{\footnotesize{
    Left: Lebesgue constants of some pointsets on the unit disk $B_2$,  for degrees $n=1,\dots,25$: AFP (red dots), DLP (blue dots), ALP (cyan dashes), Carnicer-Godes (black dashes), Halton points (magenta dashes). Right: the four lowest curves in detail. In these experiments, $m=4$ with a relative error $\approx 8.6\%$.}}
    \label{leb_disk}
\end{figure}

\begin{figure}[h!]
    \centering
        {\includegraphics[scale=0.3]{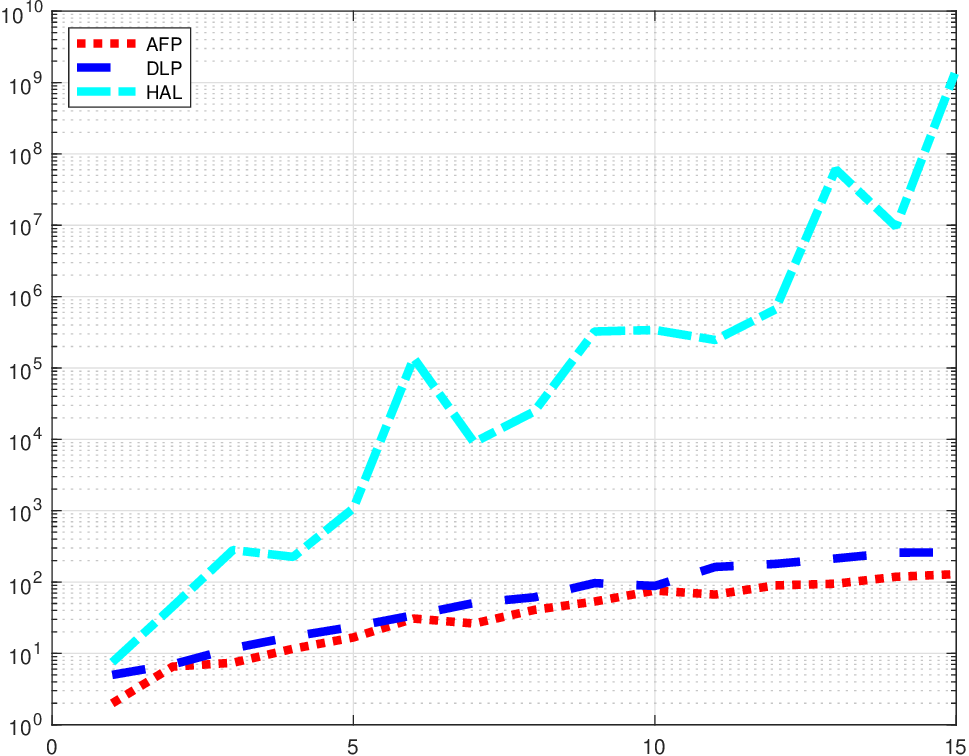}}
        {\includegraphics[scale=0.3]{leb_cube2.eps}}
    \caption{\footnotesize{
    Left: Lebesgue constants of some pointsets on the cube $[-1,1]^3$,  for degrees $n=1,\dots,15$: AFP (red dots), DLP (blue dashes), Halton points (cyan dashes).
    Right: the two lowest curves in detail. In these experiments, $m=3$ with a relative error $\approx 7.7\%$.}}
    \label{leb_cube}
\end{figure}

\begin{figure}[h!]
    \centering
        {\includegraphics[scale=0.3]{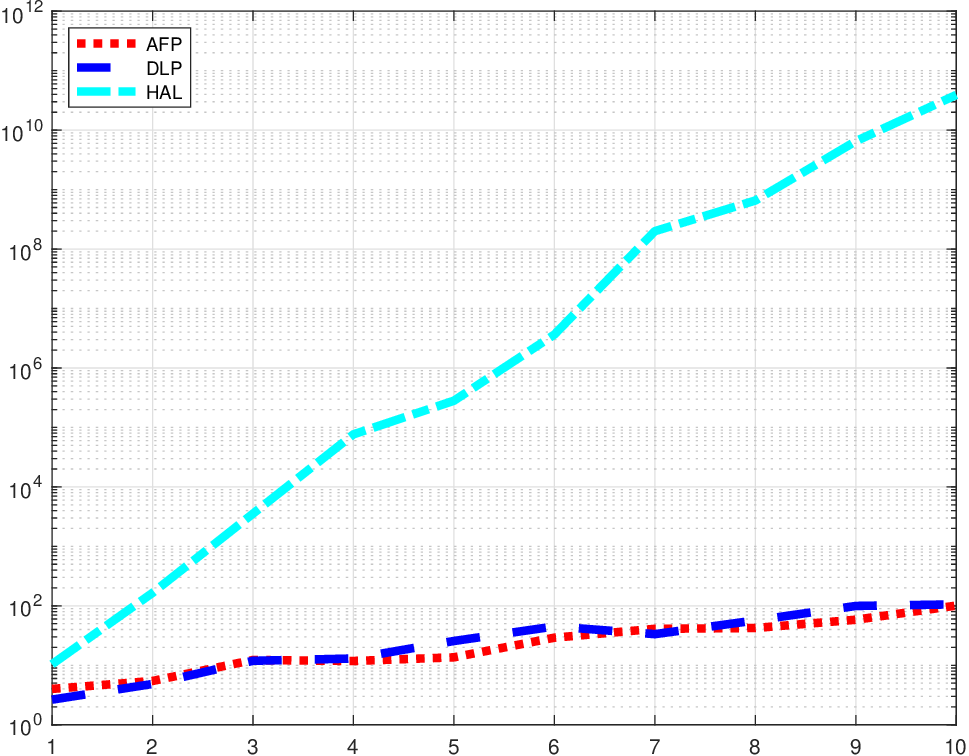}}
        {\includegraphics[scale=0.3]{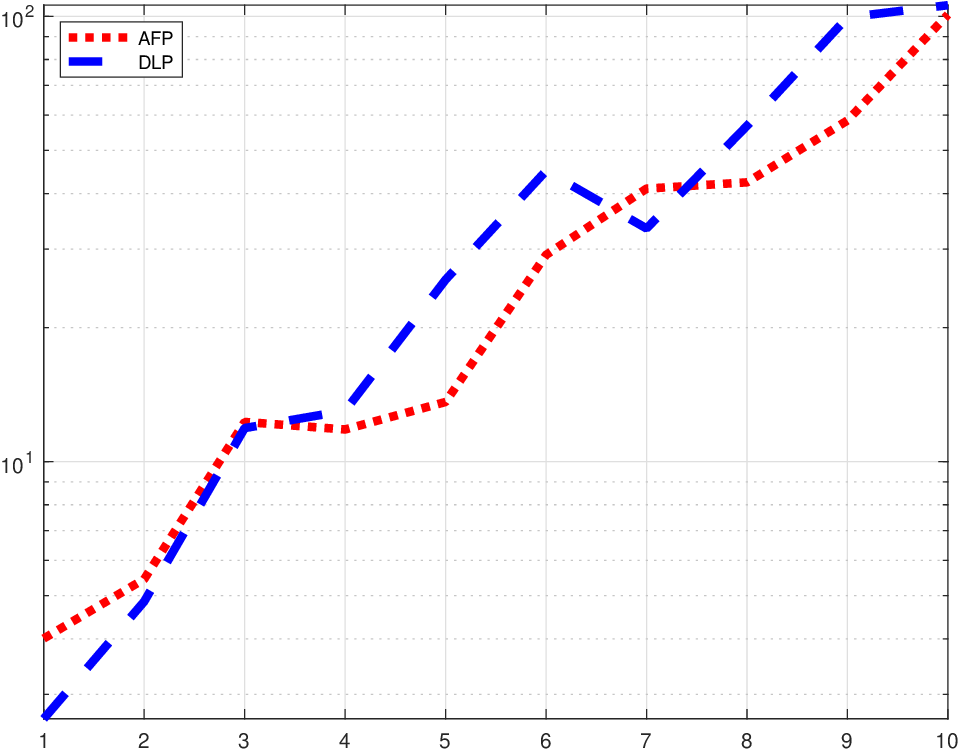}}
    \caption{\footnotesize{
    Left: Lebesgue constants of some pointsets on the unit ball $B_3$,  for degrees $n=1,\dots,10$: AFP (red dots), DLP (blue dashes), Halton points (cyan dashes).
    Right: the two lowest curves in detail. In these experiments, $m=5$ with a relative error $\approx 8.1\%$.}}
    \label{leb_ball}
\end{figure}

For the purpose of illustration, in Figures \ref{leb_disk}-\ref{leb_ball} 
we plot the computed Lebesgue constants for interpolation at AFP, DLP and Halton points on disk, cube and ball. The AFP and DLP have been computed on the same polynomial meshes $X_n^m$ used for the Lebesgue constant evaluation. On the disk we also consider the Carnicer-Godes interpolation points \cite{CG14} and the Approximate Lebesgue Points (ALP) computed in \cite{MS19}. 

\vskip0.2cm

\vskip0.2cm

\vskip0.2cm

\subsection{Standard and weighted least-squares}
As observed in Remark 2, discrete least-squares operators, one of the very basic tools of computational mathematics in both the standard (equally weighted) and the 
weighted case, fall into the class of projectors where we can evaluate the Lebesgue constant, i.e. their uniform norm, by polynomial meshes. 

It is worth recalling that a connection of (standard) least-squares with polynomial meshes was already pointed out in \cite{CL08}, where 
it  was shown that if the sampling set is a polynomial mesh on $K$, say $\Xi=X_{n}$ with constant $c$, then $\|L_n\|\leq c\sqrt{card(X_n)}$. This is however 
only a rough bound (as there observed), whereas estimating the actual size is important in applications. 

Concerning weighted least-squares, we may also recall 
that they include for example hyperinterpolation operators (see Remark 2), as well as instances coming from the recent topic of ``compression'' of discrete measures. Roughly summarizing, given a discrete measure with large support, such a compression corresponds to extract from the support a subset of re-weighted points, such that the corresponding discrete measure keeps the same polynomial moments up to a given degree.  
This topic has been receiving an increasing attention in the literature of the last decade, in both the probabilistic and the deterministic setting; cf., e.g., \cite{BPV20,ESV22,H21,LL12,PSV17,SV15,Tche15} with the references therein.

In particular, for discrete least-squares approximation of degree $n$ on a sampling set $\Xi$, moment matching has to be imposed up to degree $2n$, thus preserving orthogonal polynomials and reproducing kernels. This can be obtained by seeking a sparse nonnegative solution to the underdetermined moment-matching system $V_{2n}^t(\Xi)w=V_{2n}^t(\Xi)u$, where $u=(1,\dots,1)^t$. Such a solution with no more than $N_{2n}=dim(\mathbb{P}_{2n})$ nonzero components exists, by the well-known Caratheodory theorem on conical combinations applied to the columns of the matrix, and can be computed by solving the NonNegative Least Squares (NNLS) 
problem 
\begin{equation} \label{NNLS}
\min_{u\geq 0}\|V_{2n}^t(\Xi)w-V_{2n}^t(\Xi)u\|_2
\end{equation}
via the Lawson-Hanson NNLS-solver \cite{LH95} and its accelerated variants, such as that based on the recently developed ``deviation maximization'' criterion instead of column pivoting in the underlying $QR$ factorizations (cf. \cite{DODM23,DM22,DMV20,DMV20-2,Slawski}). The nonzero components of $w$ then determine a compressed support $\Xi_n\subset \Xi$ with $N_n\leq card(\Xi_n)\leq N_{2n}$, where the weighted least-squares polynomial can be computed; cf. \cite{PSV17}.

\begin{figure}[h!]
    \centering
        {\includegraphics[scale=0.3]{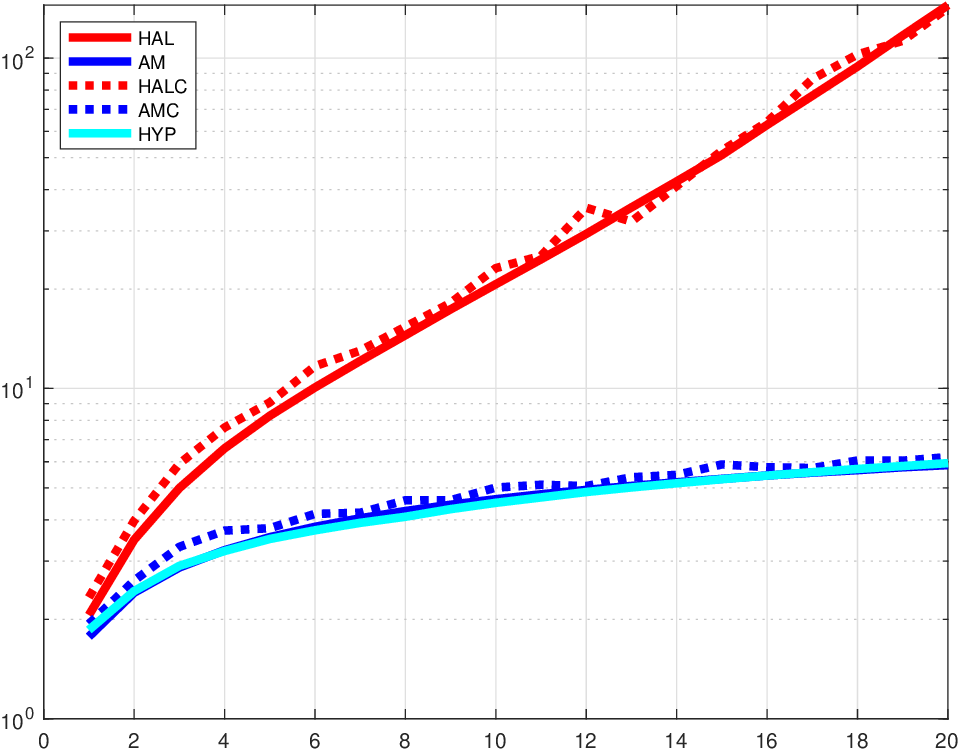}}
        {\includegraphics[scale=0.3]{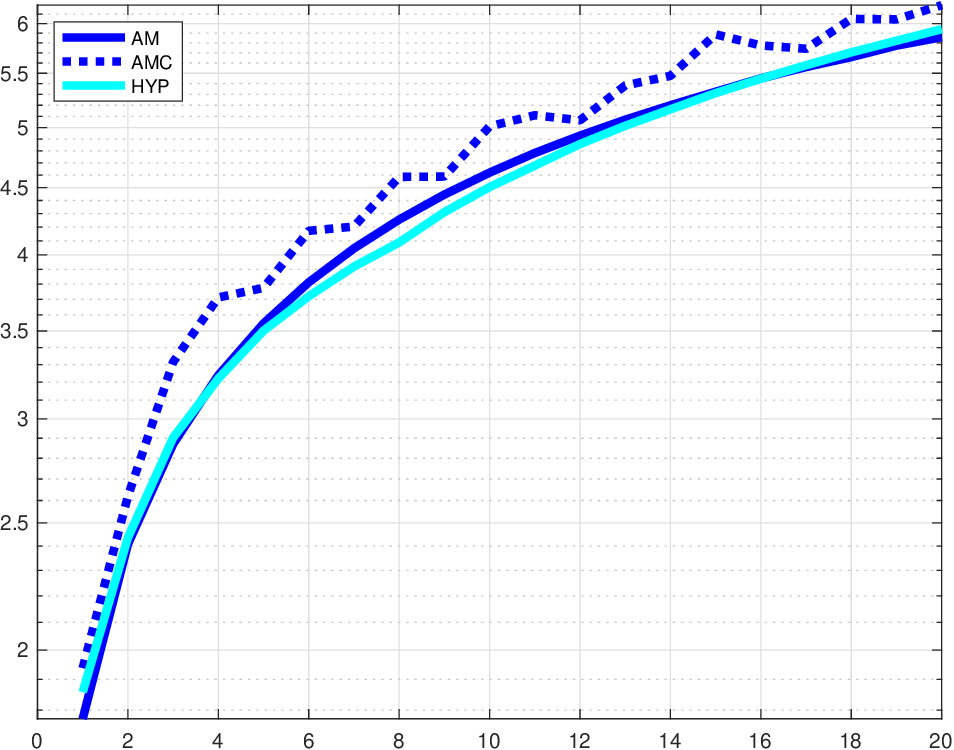}}
    \caption{\footnotesize{
    Left: Lebesgue constants of least-squares on the square $[-1,1]^2$,  for degrees $n=1,\dots,20$, on: 10000 Halton points (HAL, red dots) and compressed version (HALC, red dashes); Chebyshev mesh $X_{20}^5=\mathcal{C}_{100}\times \mathcal{C}_{100}$ with 10000 points (blue line, AM20) and  compressed version (AM20C, blue dashes); hyperinterpolation with {\color{blue}product Gauss-Chebyshev quadrature} (HYP, cyan line). Right: the three lowest curves in detail. 
    In these experiments, $m=3$ with a relative error $\approx 7.7\%$.}}
    \label{leb_LS_square}
\end{figure}

\begin{figure}[h!]
    \centering
        {\includegraphics[scale=0.3]{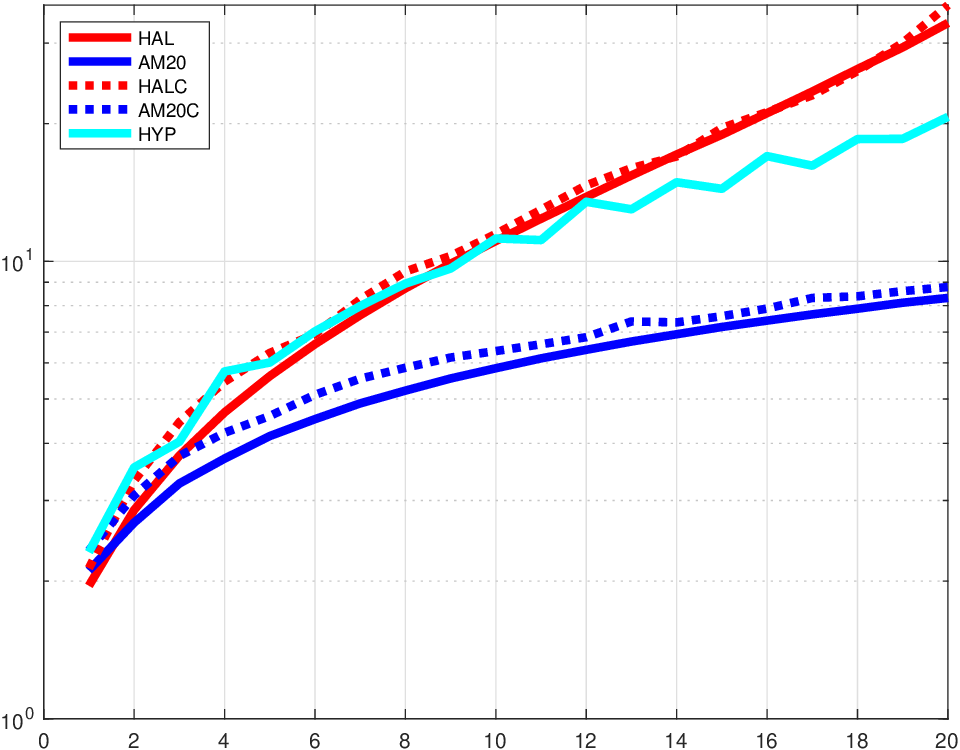}}
        {\includegraphics[scale=0.3]{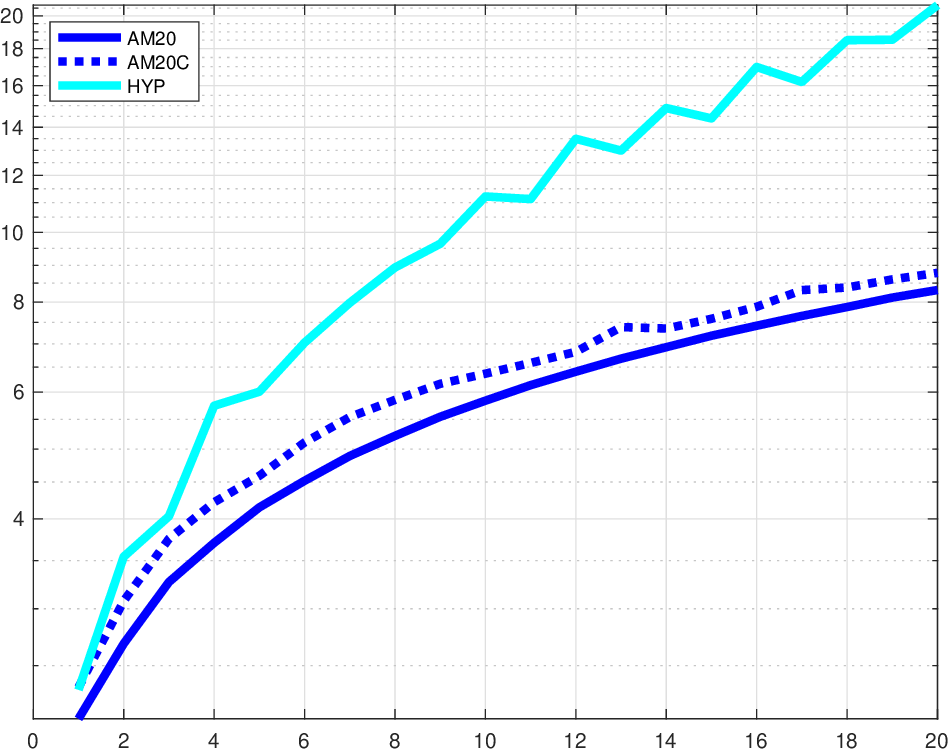}}
    \caption{\footnotesize{
    Left: Lebesgue constants of least-squares on the unit disk $B_2$, for degrees $n=1,\dots,20$, on: 12800 Halton points (HAL, red dots) and compressed version (HALC, red dashes); Chebyshev mesh $X_{20}^4=\mathcal{S} (\mathcal{C}_{80}\times \mathcal{C}_{160})$ with 12800 points (AM20, blue line) and  compressed version (AM20C, blue dashes); hyperinterpolation via low cardinality rules for the Lebesgue measure (HYP, cyan line). Right: the three lowest curves in detail.
    In these experiments, $m=4$ with a relative error $\approx 8.6\%$.}}
    \label{leb_LS_disk}
\end{figure}

For the purpose of illustration, in Figures (\ref{leb_LS_square})-(\ref{leb_LS_disk}) we compare on square and disk the Lebesgue constants of hyperinterpolation and of polynomial least-squares on Halton points and on Chebyshev polynomial meshes for degrees $n=1,\dots,20$, together with their compressed versions. To give an idea of the compression ratios, those on the disk are reported in Table 2. We recall that in the case of hyperinterpolation it is theoretically known that the Lebesgue constant is $\mathcal{O}(n^2)$ for the square and $\mathcal{O}(n)$ for the disc with the Lebesgue measure, and $\mathcal{O}(\log^2(n))$ for the square with the product Chebyshev measure; cf. e.g. \cite{DMSV14,Wade13,WWW14}. 

It is numerically manifest that Lebesgue constants of full and compressed least-squares have substantially the same size, with a remarkable reduction of the sampling cardinality for the latter. Moreover, the Lebesgue constant of least-squares on Chebyshev meshes and of hyperinterpolation with the product Chebyshev measure turn out to be very close (see Figure \ref{leb_LS_square}). This is not really surprising, since the uniform discrete measure 
supported at univariate Chebyshev points is an algebraic quadrature formula for the (normalized) Chebyshev measure, and such a behavior extends to a product-like framework. Hence, standard discrete least-squares on a Chebyshev mesh of the square are equivalent to an hyperinterpolation with respect to the product Chebyshev measure, whose Lebesgue constant is expected to be $\mathcal{O}(\log^2(n))$ in view of \cite{WWW14}.
The same can be said for the compressed least-squares, since they also correspond to an algebraic quadrature with the same moments up to degree $2n$. Notice that the Lebesgue functions of standard and compressed least-squares are not coincident, but both correspond to  hyperinterpolation with respect to the product Chebyshev measure.

A similar argument applies in interpreting Figure \ref{leb_LS_disk}, since standard discrete least-squares on a polar Chebyshev mesh of the disk are equivalent to an hyperinterpolation with respect to its equilibrium measure $dx_1\,dx_2/(\pi\sqrt{1-x_1^2-x_2^2})=r\,dr\,d\theta/(\pi\sqrt{1-r^2})$, whose Lebesgue constant is expected to be $\mathcal{O}(\sqrt{n})$ in view of \cite{Wade13}.

\vskip0.2cm

\vskip0.2cm

\begin{table}[ht]
\begin{center}
{\footnotesize
\begin{tabular}{|c|c|c|c|c|c|c|c|c|c|c|}
\hline
$n$ & 2 & 4 & 6 & 8 & 10 & 12 & 14 & 16 & 18 & 20\\
\hline \hline
$card=N_{2n}$ & 15 & 45 & 91 & 153 & 231 & 325 & 435 & 561 & 703 & 861\\
\hline
$cmp\;ratio$ & 853 & 284 & 141 & 84 & 55 & 40 & 29 & 23 & 18 & 15\\
\hline
\end{tabular}
}
\caption{\footnotesize{Cardinalities and sampling compression ratios for compressed polynomial LS of degree $n$ on 12800 points of the disk.}}
\label{relerrors}
\end{center}
\end{table}

\vskip0.8cm 
\noindent
{\bf Acknowledgements.} 
Work partially supported by the DOR funds of the University of Padova, by the INdAM-GNCS grant ``Multivariate approximation and integration with application to integral equations" (A. Sommariva, M. Vianello), and 
by the National Science Center - Poland, grant Preludium Bis 1, N. 2019/35/O/ST1/02245 (D.J. Kenne). The research cooperation was funded by the program Excellence Initiative – Research University at the Jagiellonian University in Krak\`{o}w (A. Sommariva). 
This research has been accomplished within the RITA ``Research ITalian network on Approximation" and the SIMAI Activity Group ANA\&A (A. Sommariva, M. Vianello), and the UMI Group TAA ``Approximation Theory and Applications" (A. Sommariva).

\end{document}